\UseRawInputEncoding
\documentclass[12pt,leqno]{article}
\usepackage{enumerate,amsmath,amsthm,amssymb,latexsym,
amsfonts,amscd}

\newtheorem{defn}{Definition}[section]
\newtheorem{thm}[defn]{Theorem}
\newtheorem{cor}[defn]{Corollary}
\newtheorem{prop}[defn]{Proposition}
\newtheorem{lem}[defn]{Lemma}

\def\F{\mbox{${\cal F}$}}

\def\O{\mbox{${\cal O}$}}
\newcommand{\rr}{\mathbb R}
\newcommand{\cc}{\mathbb C}
\newcommand{\zz}{\mathbb Z}
\newcommand{\nn}{\mathbb N}

\input{arrow.tex}

\begin{document}
\title{\bf On functional analytic approach for Gleason's problem in the theory of SCV}
\author{S. R. PATEL}
\maketitle
\begin{tabbing}
\> Secondary: 13F25, 32A05,\\
\> 32A17, 46E25\\
\end{tabbing}
\newpage
\noindent {\bf Abstract.}  This paper establishes the Gleason
result for finitely generated ideals in the context of Fr\'{e}chet
algebras, and, in particular, provides an affirmative answer to a
question about the Gleason result in commutative Fr\'{e}chet
algebras (Carpenter posed this question for uniform Fr\'{e}chet
algebras in 1970). As a welcome bonus of our method, {\it locally
Stein algebras} are completely characterized, and, as an
application of this characterization, an affirmative answer to the
Gleason problem (1964) for such algebras is provided recapturing
all the classical results on the Gleason problem in the
theory of several complex variables.\\
\newpage
\section{Introduction and statement of the main theorem} An important
subject in the theory of Fr\'{e}chet algebras is the question of
the existence of analytic structure in spectra. The detailed study
of this problem for uniform Fr\'{e}chet algebras is discussed in
~\cite{15}, especially the work of Brooks, Carpenter, Goldmann and
Kramm. In particular, several function algebraic characterizations
of certain types of Stein algebras are given (e.g., see
~\cite[Thm. 11.1.4]{15}; Kramm's Theorem). As far as we know, no
characterizations are obtained by studying the ideal structure of
a Fr\'{e}chet algebra since 1964. In fact, the solution of the
Gleason problem for locally Stein algebras (defined below),
obtained through the functional analytic approach, recaptures all
the classical results on the Gleason problem in the theory of SCV
(see 4.4 below).

We remark that an ideal $I$ of a Fr\'{e}chet algebra $A$ may be
finitely generated in two senses; it may be the least ideal
containing a certain finite set $F$ or it may be the least closed
ideal containing $F$. The former situation, called algebraically
finitely generated, is rare, and the latter case is common. Ideals
which are algebraically finitely generated and closed at the same
time are distinctly uncommon in Fr\'{e}chet algebras. However,
algebraically finitely generated ideals are not trivial, and do
form an analytic variety in a natural way provided that they
belong to the spectrum of $A$. Also, very rarely is it possible to
give an analytic structure to the whole spectrum and therefore it
is of interest to know conditions which ensure that parts of the
spectrum of a Fr\'{e}chet algebra can be given an analytic
structure. In ~\cite{18}, Loy gave a sufficient condition for the
existence of local analytic structure in spectra of certain
commutative Fr\'{e}chet algebras.

However, this current work is specifically concerned with the
determination of sufficient conditions for the existence of local
analytic structure in the spectrum of a commutative Fr\'{e}chet
algebra by studying the structure of the algebra. (See Main
Theorem below.) This approach was taken in ~\cite{24} for
principal ideals, and a very general result was given for locally
Riemann algebras, providing extensions of results in ~\cite{7, 19}
and elsewhere. This extends the result of Gleason to finitely
generated ideals in Fr\'{e}chet algebras, answering affirmatively
a question posed in ~\cite{7} for Fr\'{e}chet algebras. As a
consequence, locally Stein algebras are completely characterized
by intrinsic properties within the class of Fr\'{e}chet algebras
(see Thm. 4.1 below). Though the present paper is primarily
addressed to functional analysts; we hope that complex analysts
may also find the {\it sufficient conditions} interesting from
applications point of view, for example, see Theorems 4.2, 4.3 and
Corollary 4.4 below, and references (and their reviews, too) to
~\cite{14} in MathSciNet (MR0159241 (28 \#2458)). Also, it is
worthwhile mentioning some work that is somewhat (tangentially)
related to the work given here (see ~\cite{3, 9, 21, 26} for more
details).

Throughout the paper, ``algebra" will mean a non-zero, complex
commutative algebra with identity $e$. We recall that a {\it
Fr\'{e}chet algebra} is a complete, metrizable locally convex
algebra $A$ whose topology $\tau$ may be defined by an increasing
sequence $(p_m)_{m\in \nn}$ of submultiplicative seminorms. The
basic theory of Fr\'{e}chet algebras was introduced in ~\cite{15,
20}. The principal tool for studying Fr\'{e}chet algebras is the
Arens-Michael representation, in which $A$ is given by an inverse
limit of Banach algebras $A_m$. We shall briefly describe this in
Section 2, in order to establish notation that will be used
throughout the paper. A Fr\'{e}chet algebra $A$ is called a {\it
uniform Fr\'{e}chet algebra} if for each $m\,\in\,\nn$ and for
each $x \,\in\,A$, $p_m(x^2)\,=\,p_m(x)^2$. In Section 2 we define
the concept of a Fr\'{e}chet algebra of power series in $\F_k$.
The proof of the main theorem, presented in Section 3, is broken
up into several technical results of some independent interest.
The paper ends with some remarks on the hypotheses of the main
theorem.

\indent We recall that the spectrum $M(A)$ (with the Gel'fand
topology) has an {\it analytic variety} at $\phi \,\in \,M(A)$ if
there is a subvariety $D$ containing $0$ of a domain in some
$\cc^k$ and a continuous injection $f : D \rightarrow M(A)$ such
that $f(0) = \phi$ and $\hat{x} \circ f \in \textrm{Hol}(D)$ for
all $x \,\in\,A$. Also, if $\phi\,\in\,M(A)$, then we say that
$\hat{A}$ has an analytic structure in $\phi$ provided there are
open neighbourhood $U$ of $\phi$ in $M(A)$, an analytic subset $Z$
of a domain $G$ in some $\cc^k$ and a homeomorphism $p: Z
\rightarrow U$ such that $\hat{x}\circ p \in \textrm{Hol}(Z)$ for
all $x \in A$. A uniform Fr\'{e}chet algebra $A$ is called a {\it
Stein algebra} if it is topologically and algebraically isomorphic
to the Fr\'{e}chet algebra $\textrm{Hol}(X)$ of all holomorphic
functions on some (reduced) Stein space $X$ (see
~\cite[11.1.1]{15}). We call a Fr\'{e}chet algebra $A$ a {\it
locally Stein algebra} if a non-empty part of $M(A)$ can be given
the structure of a (reduced) Stein space in such a way that the
completion in the compact open topology of the algebra of Gel'fand
transforms of elements of $A$, restricted to this part, is the
Fr\'{e}chet algebra of all holomorphic functions on this Stein
space. We remark that a locally Stein algebra $A$ is not
necessarily semi-simple; e. g., locally Stein algebras
$\textrm{Hol}(X) \,\times\,B$, where $B$ is a non-semi-simple
Fr\'{e}chet algebra and $X$ is a reduced Stein space.

Let $x \in A$ and let $R_x$ denote the linear operator $A
\rightarrow A$ of multiplication by $x$. A non-zero element $x \in
A$ is called a {\it strong topological divisor of zero} in $A$ if
$R_x$ is not an isomorphism into, i.e. a linear homeomorphism of
$A$ onto $Ax$ ~\cite[Def. 11.1]{20}. A non-zero element $x \in A$
is a {\it topological divisor of zero} in $A$ if, whenever a
sequence $(p_m)$ of submultiplicative seminorms defines the
Fr\'{e}chet topology of $A$, there exists $l \in \nn$ such that
$x_l$ is a topological divisor of $0$ in $A_l$ ~\cite[Def.
11.2]{20}. The two notions of topological divisor of zero agree
for normed algebras.

We now state the main result on analytic structure to be proved in
this paper. We note that: (i) there is a one-to-one correspondence
$\phi \rightarrow \ker \phi$ between the set $M(A)$ of all
non-zero continuous complex homomorphisms on $A$ and the closed
maximal ideals in $A$, and (ii) maximal ideals of Fr\'{e}chet
algebras are not, in general, closed ~\cite[Ex., p. 83]{15};
however, every algebraically finitely generated, maximal ideal is
a closed maximal ideal in $A$ ~\cite{2}. \linebreak For each $n
\in \nn$, $M^n$ is the ideal generated by products of $n$ elements
in $M$.

\noindent {\bf MAIN THEOREM.} {\it Let $A$ be a commutative,
unital Fr\'{e}chet algebra, with its topology defined by a
sequence $(p_m)$ of norms and with the corresponding Arens-Michael
isomorphism $A\,=\,\lim\limits_{\longleftarrow} (A_m;\,d_m)$.
Suppose that $A$ has a maximal ideal $M$ that is algebraically
finitely generated, say by $t_1,\,t_2,\,\dots ,\,t_k$, that for
each $n$ the homogeneous monomials of degree $n$ in
$t_1,\,t_2,\,\dots ,\,t_k$ are representatives of a basis for
$\overline{M^{n}}/\overline{M^{n + 1}}$; and that the generators
$t_i$, $i\,=\,1,\,2,\,\dots ,\,k$, have the property that
$t_{i_m}$ is not a topological divisor of zero in $A_m$ for all
sufficiently large $m$. Then:
\begin{enumerate}
\item[{\rm (i)}] $A/\bigcap _{n\geq 1}\overline{M^{n}}$ is a
semisimple Fr\'{e}chet algebra of power series in $\F_k$;
\item[{\rm (ii)}] there is an analytic variety at $\phi$, where
$M\,=\,\ker \phi$; \item[{\rm (iii)}] for each $x\,\in \,A$,
$\hat{x}$ vanishes on a neighbourhood of $\phi$ provided that
\linebreak $x\,\in \,\bigcap _{n\geq 1}\overline{M^{n}}$.
\end{enumerate}}
\section{Fr\'{e}chet algebras of power series in $\F_k$} Let $A$ be a Fr\'{e}chet algebra, with its topology defined by an
increasing sequence $(p_{m})_{m\geq 1}$ of submultiplicative
seminorms. For each $m$, let $Q_{m}\,:\, A\, \rightarrow \,
A/\textrm{ker}\,p_{m}$ be the quotient map. Then
$A/\textrm{ker}\,p_{m}$ is naturally a normed algebra, normed by
setting $\|x\,+\,\textrm{ker}\,p_m\| _m\,=\, p_m(x)\;(x\,\in
\,A)$. We let $(A_m;\,\|\cdot\|_m)$ be the completion of
$A/\textrm{ker}\,p_m$; henceforth we consider $Q_m$ as a mapping
from $A$ into $A_m$. Then
$d_m(x\,+\,\textrm{ker}\,p_{m+1})\,=\,x+\textrm{ker}\,p_m\;(x\,\in\,A)$
extends to a norm decreasing homomorphism $d_{m}\,:\,A_{m +
1}\,\rightarrow\,A_{m}$ such that
$$
\harrowlength=20pt \varrowlength=40pt \sarrowlength=\harrowlength
\commdiag{A_1&\mapleft^{d_1}&A_2&\mapleft^{d_2}&
A_3&\mapleft&\cdots&\mapleft&A_m&\mapleft^{d_m}&A_{m+1}&\mapleft&\cdots}
$$
is an inverse limit sequence of Banach algebras; and
bicontinuously $A\,=\,\lim\limits_{\longleftarrow} (A_m;\,d_m)$.
This is called an {\it Arens-Michael representation} of $A$. For
an element $x\,\in\,A$, we may write $x_m\,=\,Q_m(x)$; it is then
evident that, for each $x\,\in\,A$, the sequence $(x_m)$ is an
element of $\lim\limits_{\longleftarrow} (A_m;\,d_m)$.

Let $M$ be a closed ideal of a Fr\'{e}chet algebra $A$. Then
$\overline{M^{n}}$ for each $n\,\geq \,1$ and $\bigcap _{n\geq
1}\overline{M^{n}}$ are also closed ideals of $A$. We now state
our two vital technical lemmas (see ~\cite[p. 127]{23}), recalling
the Arens-Michael representations of $M,\,\overline{M^{n}}$ for
each $n\,\geq \,1$ and $\bigcap _{n\geq 1}\overline{M^{n}}$, and
their quotient Fr\'{e}chet algebras $A/\overline{M^{n}}$ for each
$n\,\geq \,1$ and $A/\bigcap _{n\geq 1}\overline{M^{n}}$.
\begin{lem} \label{Lemma 2.1_SM} Let $M$ be a closed  ideal of $A$.
Then the Arens-Michael isomorphism $A\:\cong
\:\lim\limits_{\longleftarrow}(A_{m};\,d_{m})$ induces
isomorphisms:
\begin{enumerate}
\item[{\rm (i)}] $M\,\cong
\,\lim\limits_{\longleftarrow}(M_{m};\,\overline{d_{m}})$;
\item[{\rm (ii)}] $\overline{M^{n}} \,\cong
\,\lim\limits_{\longleftarrow}(\overline{M_{m}^{n}};\,\overline{d_{m}})\;\;(n\,\geq
\,1)$;
 \item[{\rm
(iii)}] $\bigcap _{n\geq 1}\overline{M^{n}}\,\cong
\,\lim\limits_{\longleftarrow}(\bigcap_{n\geq
1}\overline{M_{m}^{n}};\,\overline{d_{m}})$.
\end{enumerate}
(Here $\overline{d_{m}}\,=\,d_{m}\mid _{I_{m + 1}}\,:\,I_{m +
1}\,\rightarrow\,I_{m}$, where $I_m\,=\,\overline{Q_m(I)}$
(closure in $A_m$), whenever $I$ is a closed ideal in $A$.)
$\hfill \Box$
\end{lem}
\begin{lem} \label{Lemma 2.2_SM} With the above notation, the
Arens-Michael isomorphism $A\,\cong
\,\lim\limits_{\longleftarrow}(A_{m};\,d_{m})$ induces
isomorphisms:
\begin{enumerate}
\item[{\rm (i)}] $A/\overline{M^{n}} \,\cong
\,\lim\limits_{\longleftarrow}(A_{m}/\overline{M_{m}^{n}};\,\tilde{d_{m}})\;\;(n\,\geq\,1)$;
\item[{\rm (ii)}] $A/\bigcap _{n\geq 1}\overline{M^{n}}\,\cong
\,\lim\limits_{\longleftarrow}(A_{m}/\bigcap_{n\geq
1}\overline{M_{m}^{n}};\,\tilde{d_{m}})$.
\end{enumerate}

(Here $\tilde{d_{m}}\,:\,A_{m + 1}/\overline{M_{m +
1}^{n}}\,\rightarrow\,A_{m}/\overline{M_{m}^{n}}$ is the
homomorphism induced by $d_{m}$.) $\hfill \Box$
\end{lem}

Let $M$ be a closed maximal ideal of a Fr\'{e}chet algebra $A$. We
shall suppose from now on that
$\textrm{dim}(M/\overline{M^{2}})\,=\,k$ is finite (it is easy to
see that for finitely generated Fr\'{e}chet algebras this
condition is automatically satisfied; see ~\cite[Prop. 2.2]{26}
for the Banach case). Then, by the remark following Theorem 2.3 of
~\cite{26}, for each $n\,\in\,\nn$, the homogeneous monomials of
degree $n$ in $t_1,\,t_2,\,\dots ,\,t_k\,\in\,M$ are
representatives of a basis for $\overline{M^{n}}/\overline{M^{n +
1}}$ if and only if $\textrm{dim}(\overline{M^{n}}/\overline{M^{n
+ 1}})\,=\,C_{n+k-1,n}$ for all $n$, and thus $M$ is not
nilpotent. This situation arises for separable Fr\'{e}chet
algebras of power series in $\F_k$ (see ~\cite[Thm. 3.1]{25}),
and, in particular, for the uniform closure of the polynomials on
the (closed or open) unit poly-disc in $\cc^k$. Thus, in a special
case, we have the following
\begin{prop} \label{Proposition 2.3_SM} Let $(A,\,(p_m))$ be a
commutative, unital Fr\'{e}chet algebra with the Arens-Michael
isomorphism $A\,\cong
\,\lim\limits_{\longleftarrow}(A_{m};\,d_{m})$, and let $M$ be a
closed maximal ideal of $A$ such that: (i) $\bigcap _{n\geq
1}\overline{M^{n}} = \{0\}$ and (ii) \linebreak
$\textrm{dim}(\overline{M^{n}}/\overline{M^{n + 1}}) =
C_{n+k-1,n}$ for all $n$. Then there exist $t_1,\,t_2,\,\dots
,\,t_k\,\in\,M$ such that $\overline{M^{n}}\,=\,\overline{M^{n +
1}}\;\oplus \,\textrm{sp}\{t^I\,:\,\mid I\mid\,=\,n\}$ for each
$n\,\geq\,1$. Assume further \linebreak that each $p_m$ is a norm.
Then, for each sufficiently large $m$, $M_m$ (closure in $A_m$) is
a non-nilpotent maximal ideal of $A_m$ such that: (a) $\bigcap
_{n\geq 1}\overline{M_{m}^{n}} = \{0\}$ and (b)
$\textrm{dim}(\overline{M_m^{n}}/\overline{M_m^{n +
1}})\,=\,C_{n+k-1,n}$ for all $n$.
\end{prop}
{\bf Proof.} The first half of the proof has already been
discussed above. Assume further that each $p_m$ is a norm. First,
suppose that $M_m$ (closure in $A_m$) is nilpotent for some $m$.
Then there exists $n \, \in\,\nn$ such that $M_m^n\,=\,\{0\}$. But
then $M\,\subset\,M_m$ implies that $M$ is nilpotent, a
contradiction of the fact that $M$ is non-nilpotent. Thus it is
clear that $M_m$ is not nilpotent for each $m$, and so, by Lemma
~\ref{Lemma 2.1_SM} (ii), we have $\overline{Q_m(M)^n}\, =\,
\overline{M_m^n}\, \neq \,\{0\}$ for all $n,\,m$. Also, since
$\bigcap _{n\geq 1}\overline{M^{n}}\,=\,\{0\}$, we have $\bigcap
_{n\geq 1}\overline{M_m^{n}}\,=\,\{0\}$ for each $m$, by Lemma
~\ref{Lemma 2.1_SM} (iii) and ~\cite[Cor. A.1.25]{10}. Since $M$
is a closed maximal ideal of $A$, we have $A \,=\, M \,+\, \cc$.
Thus $Q_m(M)\, +\, \cc\, =\, M \,+\, \cc$ is dense in $A_m$, and
so also is $M_m \,+\, \cc$. Since $M_m$ is closed in $A_m$, we
have $A_m\,=\,M_m \,+\, \cc$. As it is not true that $M_m\,=\,A_m$
for infinitely many $m\,\in\,\nn$, this proves that $M_m$ is a
maximal ideal of $A_m$ for each sufficiently large $m$. Repeating
this argument for
$M\,=\,\overline{M^2}\,\oplus\,\textrm{sp}\{t_j\,:\,j\,=\,1,\,2,\,\dots
,\,k\}$ and using the fact that $\overline{Q_m(\overline{M^2})} =
\overline{M_m^2}$, we have
$\textrm{dim}\,(M_m/\overline{M_{m}^{2}})=k$ for each sufficiently
large $m$.

Since, for each sufficiently large $m$, $M_m$ is a non-nilpotent
maximal ideal of $A_m$ such that $\bigcap _{n\geq
1}\overline{M_m^{n}}\,=\,\{0\}$ and
$\textrm{dim}(M_m/\overline{M_m^2})\,=\,k$, for such $m$ we also
obtain $\textrm{dim}(\overline{M_m^{n}}/\overline{M_m^{n +
1}})\,=\,C_{n+k-1,n}$ for all $n$, by
$\textrm{dim}(\overline{M^{n}}/\overline{M^{n + 1}}) =
C_{n+k-1,n}$ for all $n$ and by the consequence of ~\cite[Thm.
1]{27}. $\hfill \Box$

We remark that, in the case where
$\textrm{dim}(M/\overline{M^{2}})\,=\,1$, one deduces \linebreak
$\textrm{dim}(\overline{M^{n}}/\overline{M^{n + 1}})\,=\,1$ for
all $n$ in ~\cite[Prop. 2.3]{23}, and so we do not require
$\textrm{dim}(\overline{M^{n}}/\overline{M^{n + 1}})\,=\,1$ for
all $n$ as a stronger hypothesis, but then we do require $M$ to be
non-nilpotent there. We have an easy counter-example (see
~\cite{25}) to show that the hypothesis that
$\textrm{dim}(\overline{M^{n}}/\overline{M^{n +
1}})\,=\,C_{n+k-1,n}$ for all $n\,\in\,\nn$ is not redundant in
the proposition above.

Let $k\,\in\,\nn$ be fixed. We write $\F_k$ for the algebra
$\cc[[X_1,X_2,\dots ,X_k]]$ of all formal power series in $k$
commuting indeterminates $X_1,X_2,\dots ,X_k$, with complex
coefficients. A fuller description of this algebra is given in
~\cite[\S 1.6]{10}, and for the algebraic theory of $\F_k$, see
~\cite[Ch. VII]{29}; we briefly recall some notations, which will
be used throughout the paper. Let $k\in \nn$, and let
$J=(j_1,j_2,\dots ,j_k)\in \zz^{+k}$. Set $\mid J\mid=j_1 + j_2 +
\cdots + j_k$; ordering and addition in \linebreak $\zz^{+k}$ will
always be component-wise. A generic element of $\F_k$ is denoted
by
$$\sum_{J\in \zz^{+k}}\lambda_J\,X^J\,=\,\sum
\{\lambda_{(j_1,\,j_2,\,\dots ,\,j_k)}\,X_1^{j_1}X_2^{j_2}\cdots
X_k^{j_k}\,:\,(j_1,\,j_2,\,\dots ,\,j_k)\,\in\,\zz^{+k}\}.$$ The
algebra $\F_k$ is a Fr\'{e}chet algebra when endowed with the weak
topology $\tau_c$ defined by the coordinate projections $\pi_I :
\F_k \rightarrow \cc, I \in \zz^{+k}$, where $\pi_I(\sum_{J\in
\zz^{+k}}\lambda_J\,X^J)=\lambda_I$. A defining sequence of
seminorms for $\F_k$ is $(p_m^{'})$, where $p_m^{'}(\sum_{J\in
\zz^{+k}}\lambda_J\,X^J)\,=\,\sum_{\mid J\mid \leq m}\mid
\lambda_J \mid\;(m\,\in\,\nn)$. A {\it Fr\'{e}chet algebra of
power series in $k$ commuting indeterminates} is a subalgebra $A$
of $\F_k$ such that $A$ is a Fr\'{e}chet algebra containing the
indeterminates $X_1,\,X_2,\,\dots ,\,X_k$ and such that the
inclusion map $A\,\hookrightarrow \,\F_k$ is continuous
(equivalently, the projections $\pi_I,\,I\,\in\,\zz^{+k}$, are
continuous linear functionals on $A$) ~\cite{11}. When it can
cause no confusion, we may use the term ``algebra of power series
in $\F_k$" for  ``algebra of power series in $k$ commuting
indeterminates" in the following; thus, Fr\'{e}chet algebras of
power series in $\F_1$ are the usual Fr\'{e}chet algebras of power
series. Though Fr\'{e}chet algebras of power series in $\F_k$
(shortly: FrAPS in $\F_k$) have been considered earlier by Loy
~\cite{18}, recently these algebras---and more generally, the
power series ideas in general Fr\'{e}chet algebras---have acquired
significance in understanding the structure of Fr\'{e}chet
algebras ~\cite{1, 6, 11, 23, 24, 25}. For examples of Fr\'{e}chet
algebras of power series, we refer to ~\cite{6}; also, analogous
examples of Fr\'{e}chet algebras of power series in $\F_k$ can be
constructed from the examples given in ~\cite{6}. In particular,
we shall consider {\it Beurling-Fr\'{e}chet algebras
$\ell^1(\zz^{+k},\,\Omega)$ of semiweight type} in the following;
we define these algebras as follows.

First, recall that $\omega$ is a {\it proper semiweight} if
$\omega(N_0)\,=\,0$ for some $N_{0}\,\in\,\nn^k$. Let
$k\,\in\,\nn$, and let
$$\ell^1(\zz^{+k},\,\Omega):=\{f=\sum_{J\in
\zz^{+k}}\lambda_J\,X^J\in \F_k:\sum_{J\in
\zz^{+k}}\mid\lambda_J\mid\,\omega_m(J)<\infty\,\textrm{for
all}\;m\},$$ where $\Omega\,=\,(\omega_m)_{m\in\nn}$ is a
separating and increasing sequence of proper semiweights on
$\zz^{+k}$ defined by $\omega_m(N)=p_m(X^N)$. Then $\rho\,=\,0$ if
and only if $\ell^1(\zz^{+k},\,\Omega)$ is a local Fr\'{e}chet
algebra if and only if the completion of
$\ell^1(\zz^{+k},\,\omega_m)/\ker p_m$ under the induced norm
$p_m$ is a local Banach algebra for all $m$. So
$\ell^1(\zz^{+k},\,\Omega)$ is either $\F_k$ or a local FrAPS in
$\F_k$. Such a Beurling-Fr\'{e}chet algebra is called
$\ell^1(\zz^{+k},\,\Omega)$ an {\it algebra of semiweight type}
(see ~\cite{25} for more details). We note that the unique maximal
ideal of $\ell^1(\zz^{+k},\,\Omega)$ is
$$\{f=\sum_{J\in \zz^{+k}}\lambda_J\,X^J\in
\ell^1(\zz^{+k},\,\Omega)\,:\,\lambda_0\,=\,0\}.$$

Let $A$ be a Fr\'{e}chet algebra of power series in $\F_k$. Then
$A$ is an integral domain. Set
$M\,=\,\textrm{ker}\,\pi_{0}\,=\,\{\sum_{J\,\in\,\zz^{+k}}\lambda_J\,X^J\,\in\,A\,:\,\lambda_0\,=\,0\}$
(where $0$ as a suffix denotes $\{0,\,0,\,\dots,\,0\}$). Then $M$
is a non-nilpotent, closed maximal ideal of $A$. Note that $\pi_0$
is a continuous projection on $A$, which is also a complex
homomorphism on $A$. Further, $\overline{M^n}\, \subset \,
\textrm{ker}\,\pi_{N - 1}$ for each $N\,\in \nn^n$, where
$M^n\,=\,\sum_{\mid N\mid\,=\,n}X^N\,A$, so that $\bigcap _{n\geq
1}\overline{M^{n}}\, =\, \{0\}$. The following result generalizes
this argument. We recall that a Fr\'{e}chet algebra of power
series in $\F_k$ satisfies {\it condition (E)} if there is a
sequence $(\gamma _K)$ of positive reals such that $(\gamma
_K^{-1}\,\pi_K)$ is an equicontinuous family ~\cite{18}, and, by
~\cite[Thm. 3.10]{25}, Fr\'{e}chet algebras of power series in
$\F_k$ (except the Beurling-Fr\'{e}chet algebras
$\ell^1(\zz^{+k},\,\Omega)$ of semi-weight type) satisfy this
condition.
\begin{thm} \label{THEOREM 2.4_P2.} Let $A$ be a Fr\'{e}chet algebra,
$\theta \,:\,A\,\rightarrow\,B$ a homomorphism of $A$ onto a
Fr\'{e}chet algebra of power series $B$ in $\F_k$ such that $B$ is
not equal to a Beurling-Fr\'{e}chet algebra
$\ell^1(\zz^{+k},\,\Omega)$ of semi-weight type. Then $A$ contains
a non-nilpotent, closed maximal ideal $M$ such that
 $\bigcap _{n\geq 1}\overline{M^{n}}\, =\, \textrm{ker}\,\theta
 $. $\hfill \Box$
\end{thm}
{\bf Remark.} We first note that the range of $\theta$ is not
one-dimensional, so, by ~\cite[Thm. 4.1]{25}, $\theta $ is
continuous whenever $B$ is not equal to the Beurling-Fr\'{e}chet
algebra $\ell^1(\zz^{+k},\,\Omega)$ of semi-weight type, and hence
one can follow the proof given in ~\cite[Thm. 1]{19} for the
Banach algebra case. In view of the further remarks on this
theorem (cf. ~\cite[REM.]{24}), it is of interest to construct
counter-examples of Fr\'{e}chet algebras which have discontinuous
epimorphisms onto the Beurling-Fr\'{e}chet algebra
$\ell^1(\zz^{+k},\,\Omega)$ of semi-weight type (and, in
particular, onto $\F_k$); cf. ~\cite[Thm. 5.5.19]{10}, and
~\cite[\S 2]{28} in which Thomas provided necessary conditions for
the existence of an epimorphism from a Fr\'{e}chet algebra onto
$\F$.

If we delete the hypothesis that $\bigcap _{n\geq
1}\overline{M^{n}} = \{0\}$ from Theorem 3.1 of ~\cite{25}, then
we have the following theorem; we shall merely sketch a proof.
\begin{thm} \label{THEOREM 2.5_P2.} Let $A$ be a commutative, unital
Fr\'{e}chet algebra. Suppose that $A$ has a closed maximal ideal
$M$ such that dim$\,(\overline{M^{n}}/\overline{M^{n +
1}})\,=\,C_{n+k-1,n}$ for each $n$. Then $A/\bigcap _{n\geq
1}\overline{M^{n}}$ is a Fr\'{e}chet algebra of power series in
$\F_k$.
\end{thm}
\noindent {\bf Proof.} Supposing $A$ satisfies the stated
condition, there exist $t_1,t_2,\dots ,t_k\in M$ such that
$\overline{M^{n}}\,=\,\overline{M^{n + 1}}\;\oplus
\,\textrm{sp}\{t^I\,:\,\mid I\mid\,=\,n\}$ for each $n\,\geq\,1$,
by Proposition 2.3. Let $x\,\in\,A$. Then a simple induction on
$n$ shows that for $n\,\geq\,1$, $$x\,=\,\sum_{\mid I\mid \leq n}
\lambda_It^I + y_n,$$ where $y_n\,\in\,\overline{M^{n+1}}$ and the
$(\lambda_I)$ are uniquely determined (in fact, by ~\cite[p.
237]{26}, $x$ has a unique partial sum of degree $n$ for each $n$
since \linebreak dim$\,(\overline{M^{n}}/\overline{M^{n +
1}})\,=\,C_{n+k-1,n}$ for all $n$). Hence the functionals
$\pi_J\,:\,x\,\mapsto \,\lambda_J$ are uniquely defined, and
linear for all $J\,\in\,\nn^k$. Thus we have a homomorphism
$\Psi\,:\,x\,\mapsto\,\sum_{I\in\zz^{+k}}\pi_{I}(x)\,t^{I}$ from
$A$ onto an algebra of formal power series in $\F_k$ with kernel
$\bigcap _{\mid I\mid \geq 0}\textrm{ker}\,\pi_{I}\,= \,\bigcap
_{n\geq 1}\overline{M^{n}}$. The inclusion $\bigcap _{\mid I\mid
\geq 0}\textrm{ker}\,\pi_{I}\,\subseteq \,\bigcap _{n\geq
1}\overline{M^{n}}$ is clear. For the reverse, suppose that
$x\,\in\,\bigcap _{n\geq 1}\overline{M^{n}}$ and that $x$ does not
belong to $\bigcap _{\mid I\mid \geq 0}\textrm{ker}\,\pi_{I}$. Let
$k$ be the least index, if there is one, such that $\pi_{\mid I
\mid = k}(x)\,\neq\,0$. Then $x\,=\,\sum_{\mid I \mid =
k}\pi_I(x)\,t^{I}\,+\,y_{k}$, where $y_{k}\,\in\,\overline{M^{k +
1}}$. So $t_{\mid I\mid\,=\,k}^{I}\,\in\,\overline{M^{k + 1}}$, a
contradiction.

For $x\,\in\,A$, let $\bar{x}$ denote the coset $x\,+\,\bigcap
_{n\geq 1}\overline{M^{n}}$. Then the mapping
$\bar{x}\,\mapsto\,\sum_{I\in\zz^{+k}}\pi_{I}(\bar{x})\,t^{I}$ is
an isomorphism from $A/\bigcap _{n\geq 1}\overline{M^{n}}$ onto an
algebra of formal power series in $\F_k$. One can now follow proof
of Theorem 3.1 of ~\cite{23}, in order to establish the theorem.
$\hfill \Box$

\indent As a corollary, we have the following result, with the
Beurling-Fr\'{e}chet algebras $\ell^1(\zz^{+k},\,\Omega)$ of
semi-weight type (including $\F_k$) as trivial examples. We note
that the polynomials in $k$ variables are dense in
$\ell^1(\zz^{+k},\,\Omega)$, and so, by ~\cite[Thm. 3.1]{25},
$M\,=\,\ker \pi_0$ is a non-nilpotent, closed maximal ideal such
that $\bigcap _{n\geq 1}\overline{M^n}\,=\,\{0\}$ and
dim$(\overline{M^n}/\overline{M^{n+1}})\,=\,C_{n+k-1,n}$ for all
$n$; the Beurling-Fr\'{e}chet algebras $\ell^1(\zz^{+k},\,\Omega)$
of semi-weight type (including $\F_k$) do not satisfy the latter
half of Proposition 2.3 since the topology of
$\ell^1(\zz^{+k},\,\Omega)$ is defined by proper seminorms.

\begin{cor} \label{COROLLARY 2.6_P2.} Let $A$ be a commutative, unital
Fr\'{e}chet algebra. Suppose that the polynomials in $e$ and
$t_1,\,t_2,\,\dots ,\,t_k$ are dense in $A$, and that\linebreak
dim$\,(\overline{M^{n}}/\overline{M^{n + 1}})$ = $C_{n+k-1,n}$ for
each $n$. Then $A/\bigcap _{n\geq 1}\overline{M^{n}}$ is a
Fr\'{e}chet algebra of power series in $\F_k$. $\hfill \Box$
\end{cor}
\section{Proof of the main theorem}
First, we prove a Banach-algebra analogue of the main theorem in
the following lemma generalizing Theorem 2.5 of ~\cite{26} (see
4.2 below). The method of proof of the lemma will be used again in
the proof of the main theorem.
\begin{lem} \label{Theorem 4_L2.} Let $A$ be a Banach algebra
which has an algebraically finitely generated maximal ideal $M$
generated by $t_1,\,t_2,\,\dots ,\,t_k$ for some $k\,\in\,\nn$.
Suppose that $\textrm{dim}(\overline{M^{n}}/\overline{M^{n +
1}})\,=\,C_{n+k-1,n}$ for all $n$. Then:
\begin{enumerate} \item[{\rm (i)}] $A/\bigcap _{n\geq 1}\overline{M^{n}}$ is
a semisimple Banach algebra of power series in $\F_k$; \item[{\rm
(ii)}] there is an analytic variety at $\phi$, where $M\,=\,\ker
\phi$; \item[{\rm (iii)}] for each $x\,\in\,A$, $\hat{x}$ vanishes
on a neighbourhood of $\phi$ provided that $x\,\in\,\bigcap_{n
\geq 1}\overline{M^n}$.
\end{enumerate}
\end{lem}
\noindent {\bf Proof.} First, we note that
$\overline{M^n}\,=\,\overline{M^{n+1}}\,\oplus\,\textrm{sp}\{t^I\,:\,\mid
I\mid\,=\,n\}$ for each $n\,\in\,\nn$, and that $\bigcap _{n\geq
1}\overline{M^{n}}$ is a closed ideal of $A$. By Theorem 2.5,
$B\,=\,A/\bigcap _{n\geq 1}\overline{M^n}$ is a Banach algebra of
power series in $\F_k$. For $x\,\in\,A$, let $\bar{x}$ denote the
coset $x+\bigcap _{n\geq 1}\overline{M^n}$ (which is, in fact,
$\sum_{I\in\zz^{+k}}\pi_I(\bar{x})t^I$). Since $B$ is an integral
domain, $\bar{t}_1,\,\bar{t}_2,\,\dots ,\,\bar{t}_k$ are certainly
not zero divisors. Also the image of $M$ under the quotient map is
an algebraically finitely generated, maximal ideal generated by
$\bar{t}_1,\,\bar{t}_2,\,\dots ,\,\bar{t}_k$, i.e., $M/\bigcap
_{n\geq 1}\overline{M^n}\,=\,\sum _{i=1}^k B\bar{t}_i$, and so, by
~\cite{2}, it is closed. Since
$\bar{t}_i,\;i\,=\,1,\,2,\,\dots,\,k,$ are in $M/\bigcap _{n\geq
1}\overline{M^n}$, they do not belong to $\textrm{Inv}\,B$, and so
each $B\bar{t}_i$ is an ideal in $B$ such that
$\bar{t}_i\,\in\,B\bar{t}_i$. In fact, each $B\bar{t}_i$ is a
closed, principal ideal in $B$ since it cannot be dense as it is
contained in a closed maximal ideal $M/\bigcap _{n\geq
1}\overline{M^n}$. Hence each $\bar{t}_i$ is not a topological
divisor of zero in $B$. Thus, for each $i$, the map
$R{_{\bar{t}_{i}}}\,:\bar{x}\, \mapsto \,\bar{x}\bar{t}_{i}$ has a
continuous inverse, and a simple induction gives \linebreak $
\parallel \pi_{I}\parallel \leq \:(2\parallel
R{_{\bar{t}_{i}}^{-1}}\parallel )^{\mid I\mid}$ for $\mid
I\mid\,\in\,{\zz^{+}}$. If $\Omega_1$ is an open poly-disc
centered at $0$ and radius $\delta <\textrm{min}_i((2\parallel
R{_{\bar{t}_{i}}^{-1}}\parallel )^{-1})$, then the map
$\bar{x}\,\mapsto\,\sum_{I\in\zz^{+k}}\pi_I(\bar{x})z^I$ of $B$
into $H(\bar{\Omega}_1)$ (the poly-disc algebra) is an algebraic
isomorphism so that $B$ must be semisimple. This proves (i).

For (ii), define functionals $\{\phi_\lambda |A :
\|\lambda\|_{\cc^k} < \delta\}$ on $A$ by $\phi_\lambda |A : x
\mapsto $ \linebreak
$\sum_{I\in\zz^{+k}}\pi_I(\bar{x_1})\lambda^I$. Then $\Gamma
:\lambda \mapsto \phi_{\lambda \delta}|A$ is easily seen to be an
analytic variety at $\phi$.

To prove (iii), let
$B{'}\,=\,A\,\oplus\,A\,\oplus\,\cdots\,\oplus\,A\,\oplus\,\cc$,
where $A$ is repeated $k$ times, and set $$\|\oplus
_{i=1}^kx_i\,\oplus\,\alpha\|\,=\,\max\{\|x_1\|,\,\|x_2\|,\dots\,\|x_k\|,\mid
\alpha\mid\}.$$ Then $B^{'}$ is a Banach algebra; the algebraic
operations are coordinatewise. For each
$\lambda\,=\,(\lambda_1,\,\lambda_2,\,\dots
,\,\lambda_k)\,\in\,\cc^k$, define linear operators on $A$
by\linebreak $T_{\lambda_i}\,:\,x\,\mapsto\,(t_i - \lambda_ie)x$,
and a linear map on $B^{'}$ by $T_\lambda\,:=\,\oplus_{i=1}^k
T_{\lambda_i}$. Since $M$ is an algebraically finitely generated,
maximal ideal of $A$ with generators $t_1,\,t_2,\,\dots ,\,t_k$,
clearly $T_0\,:\,\oplus_{i=1}^kx_i\,\mapsto\,\sum_{i=1}^kt_ix_i$
is a semi-Fredholm operator from $B^{'}$ onto $A$, with deficiency
1, and so there is $\eta\,>\,0$ such that $T_\lambda\,=\,T_0 -
\lambda e$ has deficiency $\leq\,1$ for
$\|\lambda\|_{\cc^k}\,<\,\eta$. Let
$\epsilon\,=\,\min(\delta,\,\eta)$. Then if
$\|\lambda\|_{\cc^k}\,<\,\epsilon$,
$T_\lambda(B^{'})\,\subset\,\ker \phi_\lambda$ and codim
$T_\lambda(B^{'})\,\geq\,$ codim $\ker \phi_\lambda$. So $\ker
\phi_\lambda\,=\,T_\lambda(B^{'})$.

Let $\Delta_1\,:=\,\{z\,\in\,\cc^k\,:\,\|z\|_{\cc^k}\,<\,\epsilon
\delta^{-1}\}$,
$U\,:=\,\{\psi\,\in\,M(A)\,:\,\mid\psi(t_i)\mid\,<\,\epsilon,\,i\,=\,1,\,2,\,\dots,\,k\}$.
Then from what we have just shown,
$\Gamma\,:\,\Delta_1\,\rightarrow\,U$ is a continuous bijection.
So if $x\,\in\,\bigcap_{n\geq1}\overline{M^n}$, then clearly
$\hat{x}(\psi_{\lambda
\delta})\,=\,\sum_{I\,\in\,\zz^{+k}}\pi_I(\bar{x}_1)(\lambda
\delta)^I\,=\,0$ for all $\psi_{\lambda \delta}\,\in\,U$. $\hfill
\Box$

We now give proof of the main theorem.

\noindent {\bf Proof of the main theorem.} First, we note that
$$\overline{M^{n}}\,=\,\overline{M^{n + 1}}\,\oplus\, \textrm{sp}\{t^I\,:\,\mid
I\mid\, =,n\} \,\supset\, \overline{M^{n+1}}$$ properly for each
$n\,\in\,\nn$, and so $\bigcap _{n\geq 1}\overline{M^n}$ is a
closed ideal of $A$. We have the following conclusions:

(a) $M\cong
\lim\limits_{\longleftarrow}M_{m};\overline{M^{n}}\cong
\lim\limits_{\longleftarrow}\overline{M_m^{n}};\bigcap _{n\geq
1}\overline{M^{n}}\cong \lim\limits_{\longleftarrow}\bigcap
_{n\geq 1}\overline{M_m^{n}}$ by Lemma2.1.

(b) $A/\overline{M^{n}}\,\cong
\,\lim\limits_{\longleftarrow}A_{m}/\overline{M_m^{n}};\,A/\bigcap
_{n\geq 1}\overline{M^{n}}\,\cong
\,\lim\limits_{\longleftarrow}A_{m}/\bigcap _{n\geq
1}\overline{M_m^{n}}$; $\overline{M^{n}}/\overline{M^{n + 1}}$
\linebreak $\cong
\,\lim\limits_{\longleftarrow}\overline{M_m^{n}}/\overline{M_m^{n
+ 1}}$ by Lemma 2.2. We have the last Arens-Michael isomorphism as
the ideals $\overline{M^{n}}$ are all distinct.

(c) By Theorem 2.5, $B = A/\bigcap _{n\geq 1}\overline{M^{n}}$ is
a Fr\'{e}chet algebra of power series in $\F_k$.

(d) We first recall that each $p_m$ is a norm. Then, by
Proposition 2.3, $M_{m}$ is a maximal ideal in $A_{m}$ for
sufficiently large $m$ such that
dim$(\overline{M_m^{n}}/\overline{M_m^{n + 1}})$ = $C_{n+k-1,n}$
for each $n$ since dim$(\overline{M^{n}}/\overline{M^{n + 1}})$ =
$C_{n+k-1,n}$ for each $n$. So, by Theorem 2.5, $B_{m} =
A_{m}/\bigcap _{n\geq 1}\overline{M_m^{n}}$ is a Banach algebra of
power series in $\F_k$ for sufficiently large $m$. Hence, by
passing to a suitable subsequence of $(q_m)$ defining the same
Fr\'{e}chet topology of $B$, we conclude, without loss of
generality, that each $B_m$ is a Banach algebra of power series in
$\F_k$. Thus, by ~\cite[Thm. 3.10]{25}, $B$ is not equal to a
Beurling-Fr\'{e}chet algebra $\ell^1(\zz^{+k},\,\Omega)$ of
semi-weight type, and, by ~\cite[Thm. 3.7]{25}, the topology of
$B$ is, indeed, defined by a sequence $(q_m)$ of norms. Not only
this, but, by ~\cite[Cor. 4.3]{25}, $B$ has a unique topology as a
Fr\'{e}chet algebra so that each $q_m$ can be taken as the
quotient norm induced by the norm $p_m$.

Next, following the arguments given in the proof of Lemma 3.1,
$\bar{t}_1,\bar{t}_2,\dots ,\bar{t}_k$ are certainly not zero
divisors. Also, by ~\cite{2}, $M/\bigcap _{n\geq
1}\overline{M^n}\,=\,\sum _{i=1}^k B\bar{t}_i$, is closed. Since
$\bar{t}_i,\;i\,=\,1,\,2,\,\dots,\,k,$ are in $M/\bigcap _{n\geq
1}\overline{M^n}$, they do not belong to $\textrm{Inv}\,B$ by
~\cite[Thm. 5.4]{20} (in the unital case), and so each
$B\bar{t}_i$ is an ideal in $B$ such that
$\bar{t}_i\,\in\,B\bar{t}_i$. In fact, each $B\bar{t}_i$ is a
closed, principal ideal in $B$ since it cannot be dense as it is
contained in a closed maximal ideal $M/\bigcap _{n\geq
1}\overline{M^n}$. So, for each $i$, the mapping
$R_{\bar{t}_i}\,:\,\bar{x}\,\mapsto\,\bar{x}\bar{t}_i$ of $B$ onto
$B\bar{t}_i$ is a homeomorphism, by the open mapping theorem. Now,
for each $m$ and $i$, the mapping $(R_{\bar{t}_i})_{m}
 : \bar{x} \mapsto \bar{x}\bar{t}_i$ of $(B,\,q_m)$ into
$(B\bar{t}_i,\,q_m)$ is a continuous linear transformation, being
a multiplication operator on the normed algebra $(B,\,q_m)$.
Lifting to the completions, for each $m$ and $i$,
$R{_{\bar{t}_{i_{m}}}}\,:\bar{x}_{m}\, \mapsto
\,\bar{x}_{m}\bar{t}_{i_{m}}$ of $B_m$ into $(B\bar{t}_i)_{m}$ is
continuous.

Assuming that for some $k\,\in\,\nn$ the generators
$\bar{t}_1,\bar{t}_2,\dots ,\bar{t}_k$ have the property that
$t_{i_m}$ is not a topological divisor of zero in $A_m$ for all
sufficiently large $m$, say $m\geq n$, we have $t_{i_m}$ is not a
zero divisor, and $A_mt_{i_m}$ is a closed ideal \linebreak in
$A_m$ containing $At_i$ so that $A_mt_{i_m}=(At_i)_{m}$ for all
$m\geq n$. Hence $B_m\bar{t}_{i_m}$\linebreak $=(B\bar{t}_i)_{m}$
for all $m\geq n$, and so $R{_{\bar{t}_{i_m}}}$ is surjective. In
fact, it is a homeomorphism, by the open mapping theorem. Thus,
without loss of generality, each $R{_{\bar{t}_{i_m}}}$ has
continuous inverse. In particular, by Lemma 3.1, we have
$$
\parallel \pi_{I}^{(1)}\parallel \,\leq \:(2\,\parallel
R{_{\bar{t}_{i_1}}^{-1}}\parallel )^{\mid I\mid}\ \ \
\textrm{for}\;\mid I\mid \,\in\,{\zz^{+}},
$$
where $\pi_{I}^{(1)}:B_1 \rightarrow \cc$ is the coordinate
projection on a Banach algebra of power series $B_1$ in $\F_k$.
Define $$\delta_1 := \min_i(
\parallel2R{_{\bar{t}_{i_1}}^{-1}}\parallel^{-1}) \leq
\liminf_{\mid I\mid}\parallel \pi_{I}^{(1)}\parallel ^{-1/\mid
I\mid}.$$ Now, following the arguments given in Lemma 3.1, if
$\Omega_1$ is the closed poly-disc centered at zero and radius
$\delta_1/2$, then the mapping $\vartheta_1\,:\,\bar{x}_{1}
\,\mapsto \,\sum_{i=0}^{\infty}\pi_{I}^{(1)}(\bar{x}_{1})\,z^{i}$
of $B_1$ into $\textrm{Hol}(\Omega_1)$, the standard poly-disc
algebra, is an algebraic isomorphism, and continuous by the closed
graph theorem. Since $\textrm{Hol}(\Omega_1)$ is semisimple, the
same holds for $B_1$. Clearly, the mapping
$Q_1\,:\,\bar{x}\,\mapsto\,\bar{x}_1$ of $B$ into $B_1$ is also a
continuous, injective homomorphism. This shows that $B$ is a
semisimple Fr\'{e}chet algebra of power series in $\F_k$. This
proves (i).

To prove (ii), one can follow the same arguments given in (ii) of
~\cite{24} by noticing that $\mid \lambda\mid$ should be replaced
by $\|\lambda\|_{\cc^k}$, where
$\lambda\,=\,(\lambda_1,\,\lambda_2,\,\dots,\,\lambda_k)$.

Also, to prove (iii), one can follow the same arguments given in
(iii) of Lemma 3.1. $\hfill \Box$
\section{Characterization of locally Stein algebras with an application to the Gleason problem and open questions} In 4.1 -
4.6 of ~\cite{24} (which deals with the case $k\,=\,1$), one
remarks on the hypotheses of the main theorem, with
counter-examples, showing that the assumptions, considered on
$\phi,\;t$ and $(p_m)$, cannot be dropped; nevertheless, as far as
they go, these remarks also support the case of several-variable
by considering the several-variable analogues of those
counter-examples. In particular, if $\phi$ were isolated, then
(ii) (of the main theorem and of Lemma 3.1) is clearly impossible.
Also, the remarks regarding: (i) a weaker hypothesis on
$\bar{t}_i$ in order to obtain a stronger form of the main theorem
(cf. 4.4 of ~\cite{24}), (ii) independency of the choice of
$(p_m)$ (cf. 4.5 of ~\cite{24}), and (iii) whether the sufficient
condition on the generators $t_i$ is necessary for the existence
of analytic variety at $\phi$ (cf. 4.6 of ~\cite{24}), are of
great interest (see Remark of Theorem ~\ref{corollary
4.1_necessary} below). We, here, deal with a different set of
remarks, specifically pertaining to the several-variable case.

First, we establish the following theorem generalizing Kramm's
Theorem ~\cite[Thm. 11.1.4]{15}. Since the method of proof is used
in the proof of Theorem 4.2 below, we present the proof here for
reader's convenience.
\begin{thm} \label{Theorem 11.1.4_locally Stein} Let $A$ be a Fr\'{e}chet
algebra. Then the following statements are equivalent.
\begin{enumerate} \item[{\rm (i)}] $A$ is a
locally Stein algebra; \item[{\rm (ii)}] $\hat{A}$ has an analytic
structure in each point $\phi \in Y \subset M(A)$; \item[{\rm
(iii)}] There exists an open cover $(U_i)_{i\in I}$ of $Y$
consisting of hemicompact, relatively compact and
$\hat{A}_Y$-convex subsets such that each algebra $\hat{A}_{U_i}$
is a Stein algebra.
\end{enumerate}
\end{thm}
{\bf Proof.} (i) $\Rightarrow$ (iii). Let $A$ be a locally Stein
algebra. Then $\hat{A}_Y$ is topologically and algebraically
isomorphic to $\textrm{Hol}(Y)$, $Y$ a (reduced) Stein space, so
let $T: \hat{A}_Y \rightarrow \textrm{Hol}(Y)$ be an isomorphism.
Let $\phi \in Y$ be an arbitrary point. Recall from
~\cite[2.3.2]{15} that $Y = M(\textrm{Hol}(Y))$, more precisely
the structure map $j : Y \rightarrow M(\textrm{Hol}(Y))$ defines a
homeomorphism. So the adjoint spectral map $T^* : Y \rightarrow
M(\hat{A}_Y)$ is a homeomorphism. Note in particular that $Y$ is
locally compact. Moreover, $Y$ is hemicompact since $M(A)$ is, and
$Y$ is a $\hat{A}$-convex subset of $M(A)$ by ~\cite[Rem. 4.3.4
(i)]{15}. Then, by ~\cite[Rem. 4.3.4 (iii)]{15}, $M(\hat{A}_Y)$
can be identified with $Y$ as a topological space. Since $Y$ is
locally compact, there exits an open, relatively compact,
hemicompact and $\hat{A}_Y$-convex neighbourhood $U \subset Y$ of
$\phi$ (we remark that $U$ being $\hat{A}$-convex will also do
since $U \subset Y$ is $\hat{A}_Y$-convex if and only if it is
$\hat{A}$-convex by ~\cite[4.3.7]{15}). Then $(T^*)^{-1}(U)$ is an
open and $\textrm{Hol}(Y)$-convex subset of $Y$, hence a (reduced)
Stein space. The map $$(\hat{A}_U, U) \rightarrow
(\textrm{Hol}(Y)_{(T^*)^{-1}(U)}, (T^*)^{-1}(U)), \; f \rightarrow
f \circ T^*,$$ is well-defined and defines a topological and
algebraical isomorphism. By ~\cite[2.3.6]{15} and since
$\textrm{Hol}((T^*)^{-1}(U))$ is complete, we have
$\textrm{Hol}(Y)_{(T^*)^{-1}(U)} = \textrm{Hol}((T^*)^{-1}(U))$.

(iii) $\Rightarrow$ (ii). Let $\phi \in Y \subset M(A)$ be an
arbitrary point, and let $\phi \in U_i$. By definition, there
exits a Stein space $X_i$ and a topological algebra homomorphism
$T_i : \hat{A}_{U_i} \rightarrow \textrm{Hol}(X_i)$. We can
identify $M(\textrm{Hol}(X_i))$ with $X_i$ by ~\cite[2.3.2]{15}.
By ~\cite[Rem. 4.3.4 (iii)]{15}, we can identify
$M(\hat{A}_{U_i})$ with $U_i$ as topological spaces. Thus $T_i^* :
X_i \rightarrow U_i$ defines a homeomorphism. By definition of a
reduced complex space, there exists an analytic subset $Z$ in a
domain $G$ of some $\cc^k$, an open neighbourhood $V$ of
$(T_i^*)^{-1}(\phi)$ in $X_i$ and a homeomorphism $p : Z
\rightarrow V$ so that $f \circ p \in \textrm{Hol}(Z)$ for all $f
\in \textrm{Hol}(X_i)$. Then $T_i^* \circ p : Z \rightarrow
T_i^*(V)$ is a homeomorphism onto the open neighbourhood
$T_i^*(V)$ of $\phi$. Let $\hat{g} \in \hat{A}$ be an arbitrary
element. Then $\hat{g}|_{U_i} \in \hat{A}_{U_i}$ and
$\hat{g}|_{U_i} \circ T_i^* \in \textrm{Hol}(X_i)$ by ~\cite[Rem.
11.1.1]{15}. Thus $\hat{g} \circ T_i^* \circ p \in
\textrm{Hol}(Z)$, i.e., $\hat{A}$ has an analytic structure in
$\phi$.

(ii) $\Rightarrow$ (i). For each $\phi \in Y$ choose an analytic
subset $Z_\phi$ in a domain of some $\cc^{k_\phi}$, an open
neighbourhood $U_\phi$ of $\phi$ in $Y$ and a homeomorphism
$p_\phi : Z_\phi \rightarrow U_\phi$ so that $\hat{f} \circ p_\phi
\in \textrm{Hol}(Z_\phi)$ for all $\hat{f} \in \hat{A}$. By
~\cite[Rem. 11.1.2 and 4.3.7]{15}, we can assume that each
$U_\phi$ is a hemicompact and $\hat{A}_Y$-convex set. Endow
$U_\phi$ with the complex structure so that $p_\phi$ becomes a
biholomorphic mapping, i.e., $$\textrm{Hol}(U_\phi) = \{g : U_\phi
\rightarrow \cc : g \circ p_\phi \in \textrm{Hol}(Z_\phi)\}.$$
Then $\hat{A}|_{U_\phi} \subset \textrm{Hol}(U_\phi)$, thus
$\hat{A}_{U_\phi} \subset \textrm{Hol}(U_\phi)$, since
$\textrm{Hol}(U_\phi)$ is complete with respect to the compact
open topology. $U_\phi$ is $\hat{A}_{U_\phi}$-morphically convex,
since it is $\hat{A}$-convex and by ~\cite[4.3.3]{15}, so
$(\hat{A}_{U_\phi}, U_\phi)$ satisfies the hypothesis of Rossi's
theorem ~\cite[Thm. 11.1.3]{15}. Hence $U_\phi =
M(\hat{A}_{U_\phi})$ can be given the structure of a (reduced)
Stein space so that $\hat{A}_{U_\phi} = \textrm{Hol}(U_\phi, O)$,
where $\textrm{Hol}(U_\phi, O)$ denotes the algebra of all
holomorphic functions with respect to this new structure. It is
easy to see that these new structures coincide on their
intersection by ~\cite[p. 200-201]{15}. Hence we have equipped $Y$
with the structure of a (reduced) complex space so that
$\hat{A}|_Y \subset \textrm{Hol}(Y)$. Since $Y$ is
$\hat{A}$-morphically convex, we get our result from Rossi's
theorem ~\cite[Thm. 11.1.3]{15}, i.e., $\hat{A}_Y$ is
topologically and algebraically isomorphic to $\textrm{Hol}(Y)$.
$\hfill \Box$

{\bf 4.1.} As a special case of the main theorem and Theorem
~\ref{Theorem 11.1.4_locally Stein}, locally Stein algebras are
completely characterized in the following theorems. These also
include a copy of proof of the result for the case of principal
ideals, Riemann surface and locally Riemann algebras ~\cite[Cor.
4.1]{24}. In fact, the proof of Corollary 4.1 of ~\cite{24} was
omitted there, but it turns out to be highly non-trivial for
finitely generated ideals exhibiting a significant difference
between the one-variable case and the several-variable case, as we
shall see now.
\begin{thm} \label{COROLLARY 4.1_P2.} Let $A$ be a Fr\'{e}chet
algebra, with its topology defined by a sequence $(p_m)$ of norms
and with the Arens-Michael isomorphism $A =$ \linebreak
$\lim\limits_{\longleftarrow} (A_m; d_m)$. Suppose that $M(A)$
contains an open subspace $Y$ such that every closed maximal ideal
$M$ in $Y$ is finitely generated such that: \linebreak (a)
$\textrm{dim}(\overline{M^{n}}/\overline{M^{n +
1}})\,=\,C_{n+k-1,n}$ for all $n$; and (b) the generators
$t_i,\;i\,=\,1,\,2,\,\dots ,\,k$, have the property that $t_{i_m}$
is not a topological divisor of zero in $A_m$ for all sufficiently
large $m$. Then $Y$ can be given the structure of a (reduced)
Stein space in such a way that, for each $x\,\in\,A$, the
restriction of $\hat{x}$ to $Y$ is analytic. In particular, if $Y$
is locally compact and connected, such that the above conditions
hold, then $A$ is a locally Stein algebra (that is, the completion
of $\hat{A} | Y$ with respect to the compact open topology is
topologically and algebraically isomorphic to $\textrm{Hol}(Y)$).
\end{thm}
{\bf Proof.} Let $A$, $Y \subset M(A)$ and a closed maximal ideal
$M$ in $Y$ satisfy the stated conditions. Then, by the main
theorem, there is an analytic variety at $\phi$, where $M = \ker
\phi$, i.e., there is a continuous injection $f_\phi : D_\phi
\rightarrow M(A)$ such that $f_\phi(0) = \phi$ and $\hat{x} \circ
f_\phi \in \textrm{Hol}(D_\phi)$ for all $x \in A$. Since $Y$ is
open in $M(A)$, there exists an open neighbourhood $U_\phi$ of
$\phi$ in $Y$, and so $f^{-1}_\phi(U_\phi)$ is open in $D_\phi$
containing $0$. Since $D_\phi$ is a subvariety, it is locally
compact, and so $f^{-1}_\phi(U_\phi)$ is locally compact. Hence
there exists an open neighbourhood $V_\phi$ of $0$ such that
$\overline{V}_\phi$ is compact and $\overline{V}_\phi \,\subset
\,f^{-1}_\phi(U_\phi)$; we may consider $\overline{V}_\phi$ as a
one-point compactification of $V_\phi$ as well as another
subvariety in $D_\phi$. Then the continuous injection
$f_\phi|_{\overline{V}_{\phi}} : \overline{V}_\phi \rightarrow
f_\phi(\overline{V}_\phi)$ is, in fact, a homeomorphism and
$f_\phi(\overline{V}_\phi)\, \subset\,U_\phi$. So
$f_\phi|_{\overline{V}_{\phi}} : \overline{V}_\phi \rightarrow
f_\phi(V_\phi)$ is also a homeomorphism. Thus, without loss of
generality, we can assume that $f_\phi : D_\phi \rightarrow
f_\phi(D_\phi)$ is a homeomorphism. Now follow the proof of (ii)
$\Rightarrow$ (i) of Theorem ~\ref{Theorem 11.1.4_locally Stein}.
$\hfill \Box$

In the converse direction, we have the following
\begin{thm} \label{corollary 4.1_necessary}  Let $A$ be a semi-simple locally Stein algebra.
Then every closed maximal ideal $M$ (corresponding to a point in a
(reduced) Stein space $Y$) is algebraically finitely generated
provided that $Y$ is dense in $M(A)$.
\end{thm}
{\bf Proof.} Let $A$ be a locally Stein algebra. Then $\hat{A}_Y$
is topologically and algebraically isomorphic to
$\textrm{Hol}(Y)$, $Y$ a (reduced) Stein space. Here we may assume
the Stein space $Y$ to be $(Y, \Lambda)$ by ~\cite[Thm. 2]{13}.
Let $M$ be a closed maximal ideal corresponding to a point, say
$\phi$, in a (reduced) Stein space $Y$, i.e., $M_\phi = \{x \in A
: \hat{x}(\phi) = 0\}$. Then its image $\hat{M}_\phi = \{\hat{x}
\in \hat{A} : \hat{x}(\phi) = 0\}$ is also a maximal ideal in
$\hat{A}$. So $\overline{\hat{M}_{\phi}|_Y}$ is a closed maximal
ideal in $\hat{A}_Y$. Since $\hat{A}_Y$ is topologically and
algebraically isomorphic to $\textrm{Hol}(Y, \Lambda)$, $(Y,
\Lambda)$ a (reduced) Stein space, the image of
$\overline{\hat{M}_{\phi}|_Y}$ (denoted by
$\overline{\hat{M}_{\phi}|_Y}$ again) in $\textrm{Hol}(Y,
\Lambda)$ is a closed maximal ideal. By ~\cite[Thm. 2]{13}, the
maximal ideal $(\overline{\hat{M}_{\phi}|_Y})_\phi$ of
$\Lambda_\phi$ is finitely generated (because $\Lambda_\phi$ is
noetherian) and by ~\cite[Thm. A]{8}, there exists elements $f_1,
f_2, \dots, f_k \in \textrm{Hol}(Y, \Lambda)$, the germs of which
generate $(\overline{\hat{M}_{\phi}|_Y})_\phi$ over
$\Lambda_\phi$. Then, following the proof of ~\cite[Thm. 2]{13},
we see that the functions $\hat{f}_0, \hat{f}_1, \dots, \hat{f}_k$
generate the ideal $\overline{\hat{M}_{\phi}|_Y}$ in $\hat{A}_Y$,
i.e., if $\hat{h} \in \hat{A}_Y$, then $\hat{h} = \sum_{j = 0}^{k}
a_j \hat{f}_j$ on $Y$ with coefficients in $\hat{A}_Y$. Since
$\hat{M}_{\phi}|_Y\;\subset\;\overline{\hat{M}_{\phi}|_Y}$,
$\hat{M}_{\phi}|_Y$ is also algebraically finitely generated by
the functions $\hat{f}_{0}|_Y, \hat{f}_{1}|_Y, \dots,
\hat{f}_{k}|_Y$ in $\hat{A}|_Y$. Since $Y$ is dense in $M(A)$,
$\hat{M}_\phi$ is finitely generated by the functions
$\hat{f}_{0}, \hat{f}_{1}, \dots, \hat{f}_{k}$ in $\hat{A}$. Since
$A$ is a semi-simple locally Stein algebra, $M_\phi$ is finitely
generated by $f_0, f_1, \dots, f_k$ in $A$, i. e., if $x \in
M_\phi$, then $x = \sum_{j = 0}^{k} \alpha_j f_j$ with
coefficients $\alpha_j \in A$. $\hfill \Box$
\\{\bf Remark.} From Theorem ~\ref{corollary 4.1_necessary}, it is
clear that the characterization obtained in Theorems
~\ref{COROLLARY 4.1_P2.} and ~\ref{corollary 4.1_necessary}, is
not complete in the sense that the sufficient condition on the
generators $t_i$ is not a necessary condition for the existence of
analytic variety at $\phi \in Y$. The algebra $A^\infty(\Gamma)$
is a counter-example in the one-variable case (see 4.2 of
~\cite{24} for more details); also, as mentioned in the beginning
of this section, we may generate a suitable counter-example
$A^\infty(\Gamma^k)$ to support the case of several-variable. We
also remark that the algebra $A^\infty(\Gamma)$ exhibits a
significant difference between Stein algebras and locally Stein
algebras, since if $R$ is a compact Riemann surface, then one
obtains the trivial Riemann algebra $\cc$ (see 4.4 of ~\cite{24}).

We note that Banach algebras satisfying Lemma 3.1 (in particular,
the poly-disc algebra $A(D^k)$ and the algebra $H^\infty(U)$ of
all bounded analytic functions on some bounded domain $U$ of
$\cc^k$, which are not nuclear), $A^\infty(\Gamma^k)$ (which is
nuclear) and Stein algebras are locally Stein algebras whereas
$C^\infty(\rr^k)$ is not a locally Stein algebra. Moreover, $A =
\{f \in \textrm{Hol}(\cc) : f(0) = f(i)\;\textrm{for
all}\;i\,=\,1,\,2,\,\dots \}$ is a uniform Fr\'{e}chet algebra,
being a closed subalgebra of the Stein algebra
$\textrm{Hol}(\cc)$; it is shown in 11.1.5 of ~\cite{15} that $A$
is not Stein, but it is, indeed, locally Stein. Further, a locally
Stein algebra $\textrm{Hol}(Y)\,\times\,C^\infty(\rr)$, where $Y$
a reduced Stein space, has a closed subalgebra $C^\infty(\rr)$
which is not locally Stein. Hence it is of interest to get a
criterion when a quotient (resp., a closed subalgebra) of a
locally Stein algebra is itself locally Stein. In the literature,
there are examples of complex function algebras with no analytic
structure in their spectra, but in the case $k\,=\,2$ (see
~\cite[Rem., p. 235]{15} for more references).

{\bf 4.2.} It is easy to see that if the hypothesis (2) of Theorem
2.5 of ~\cite{26} holds, then
$\textrm{dim}(\overline{M^n}/\overline{M^{n+1}})\,=\,C_{n+k-1, n}$
for all $n$, by the remarks following Theorem 2.3 of ~\cite{26};
the Beurling-Fr\'{e}chet algebras $\ell^1(\zz^{+k},\,\Omega)$ of
semi-weight type show that the converse is not true. Thus, in a
special case, Lemma 3.1 is a generalization of Read's Theorem. In
fact, the latter part of the hypothesis (2) of Theorem 2.5 of
~\cite{26} is a topological assumption whereas the hypotheses of
Lemma 3.1 are of purely algebraic nature. In addition, by the
Hilbert-Serre theorem (see ~\cite[p. 232]{28}), there is a
polynomial $P$ of degree exactly $k - 1$ such that $P(n)\,=\,
\textrm{dim}(\overline{M^n}/\overline{M^{n+1}})\,=\,C_{n+k-1, n}$
for all large $n$, and $k$ is the dimension at the origin of the
variety $\phi$.

{\bf 4.3.} In Theorem 4.1 of ~\cite{26}, Read considered a
strictly weaker hypothesis of assuming dim$(M/M^2)$ is finite
instead a maximal ideal $M$ being algebraically finitely
generated, giving of the generalization of the Gleason's result in
the Banach algebra case (~\cite[Ex. 5.1]{26} shows that varieties
thus obtained need not be, in general, open in the Gel'fand
topology, and ~\cite[Ex. 5.3]{26} shows that such sufficient
conditions are far from being necessary, even for open analytic
structure; see 4.4 below). One naturally conjectures that the main
theorem and Theorems 4.2 and 4.3 also hold true with this strictly
weaker hypothesis.

{\bf 4.4.} To see an application of Theorem ~\ref{corollary
4.1_necessary}, we recall a problem which in the literature is
known as the Gleason problem ~\cite{14}: to decide whether the
maximal ideal in $A$, where $A$ is either $H^\infty(\Omega')$ or
$A(\Omega')$, $\Omega'$ a bounded domain in $\cc^k$, consisting of
functions vanishing at the origin is algebraically finitely
generated by the coordinate functions (see ~\cite{4} and other
references therein for more details). Further, we recall that a
domain $\Omega'$ has the Gleason $A$-property if the problem has
an affirmative solution at all points of $\Omega'$. It is
meaningful to pose this problem for locally Stein algebras in an
appropriate manner: to decide whether the closed maximal ideal
(corresponding to a point $p\,\in\,X\,\subset\,M(A)$) in $A$,
where $A$ is a locally Stein algebra, consisting of the elements
of $A$ whose Gel'fand transforms vanish at a point
$p\,\in\,X\,\subset\,M(A)$, is algebraically finitely generated.
This was actually the problem posed by Gleason in ~\cite{14}, and
he mentioned that if the maximal ideal corresponding to the origin
is algebraically finitely generated, then it is finitely generated
by the coordinate functions by Theorem 2.2 (which would obviously
hold true in the Fr\'{e}chet algebra case); see ~\cite[p.
131-132]{14} for more details. As a corollary of Theorem
~\ref{corollary 4.1_necessary}, we have the following result,
answering affirmatively the Gleason problem for locally Stein
algebras.
\begin{cor} \label{corollary_locally Stein} Let $A$ be a semi-simple locally
Stein algebra, and let $M$ be a closed maximal ideal
(corresponding to a point $\phi\,\in\,Y\,\subset\,M(A)$) in $A$
consisting of the elements of $A$ whose Gel'fand transforms vanish
at $\phi\,\in\,Y\,\subset\,M(A)$. Then $M$ is algebraically
finitely generated provided that $Y$ is dense in $M(A)$. $\hfill
\Box$
\end{cor}

We say that a subspace $X$ of $M(A)$ has the Gleason $A$-property
if the problem has an affirmative solution at all points of $X$.
The above corollary says that the dense subspace $Y$ has the
Gleason $A$-property, where $A$ is a semi-simple locally Stein
algebra; in particular, $A$ can be either $H^\infty(\Omega')$ or
$A(\Omega')$, $\Omega'$ a bounded domain in $\cc^k$ containing the
origin. Thus we have given an abstract touch to the affirmative
solution of the Gleason problem, and so the abstract method given
here recaptures the classical results obtained by Beatrous Jr.
~\cite{5}, Forn\ae ss and \O vrelid ~\cite{12}, Kerzman and Nagel
~\cite{16}, Lieb ~\cite{17}, Noell ~\cite{22} and Backlund and
F\.{\.a}llstr\.{\.o}m (see ~\cite{4} and other references therein
for a list of papers on the Gleason problem), that is, a bounded
domain $\Omega'$ in $\cc^k$ has the Gleason $A$-property, where
$A$ is either $H^\infty(\Omega')$ or $A(\Omega')$ and $\Omega'$ a
(strictly or weakly) pseudoconvex domain in $\cc^k$ with various
boundary conditions.

Next, ~\cite{12} is a good reference; in the final paragraph of \S
1, the authors state the main theorem can still be proved by
replacing $A(\Omega')$ by various holomorphic H\.older- and
Lipschitz-spaces and by replacing the coordinate functions as the
generators by arbitrary generators of the maximal ideal in these
spaces. They mention a future paper, but, as far as we know, it
has never been published. In this connection, we can consider the
above application of Theorem ~\ref{corollary 4.1_necessary}, which
establishes the similar claim, but for semi-simple locally Stein
algebras.

{\bf 4.5} We say that $A$ satisfies the weak identity theorem in
$\phi\,\in\,M(A)$ if there is a (fixed) neighbourhood $U$ of
$\phi$ in $M(A)$ such that $\hat{x}|_U\,\equiv\,0$ for each
$\hat{x}\,\in\,\hat{A}$ which vanishes in an arbitrary
neighbourhood of $\phi$. It follows from the theory of several
complex variables that $A$ satisfies the weak identity theorem in
$\phi$ if there is an analytic variety at $\phi$. Thus the
statement (iii) of the main theorem (and Lemma 3.1, too) shows
that the weak identity theorem holds for locally Stein algebras.

\noindent Address: Ahmedabad, Gujarat, INDIA, E-mails: srpatel.math@gmail.com, coolpatel1@yahoo.com\\

\begin{thebibliography}{Dillo 83}
\bibitem[1] {1} G. R. Allan, Fr\'{e}chet algebras and formal power series, Studia Math. 119 (1996), 271-278.

\bibitem[2]{2} R. Arens, Dense inverse limit rings, Michigan Math. J. 5 (1958), 169-182.

\bibitem[3]{3} R. M. Aron, B. J. Cole and T. W. Gamelin, Spectra of algebras of analytic functions on a Banach space, J. Reine Angew. Math. 415 (1991), 51-93.

\bibitem[4]{4} U. Backlund and A. F\.{\.a}llstr\.{\.o}m, The Gleason property for Reinhardt domains, Math. Ann. 308 (1997), 85-91.

\bibitem[5]{5} F. Beatrous Jr, H\.{\.o}lder estimates for the $\bar{\partial}$-equation with a support condition, Pacific J. Math. 90 (1980), 249-257.

\bibitem[6] {6} S. J. Bhatt and S. R. Patel, On Fr\'{e}chet algebras of power series, Bull. Austral. Math. Soc. 66 (2002), 135-148.

\bibitem[7]{7} R. L. Carpenter, Principal ideals in F-algebras, Pacific J. Math. 35 (1970), 559-563.

\bibitem[8]{8} H. Cartan, Vari\'{e}t\'{e}s analytiques complexes et cohomologie, Coll. de Bruxelles (1953), 41-55.

\bibitem[9]{9} B. J. Cole, T. W. Gamelin, W. B. Johnson, Analytic disks in fibers over the unit ball of a Banach space, Michigan Math. J. 39 (1992), 551-569.

\bibitem[10]{10} H. G. Dales, Banach algebras and automatic continuity, London Math. Soc. Monogr. 24 (Clarendon Press, 2000).

\bibitem[11]{11} H. G. Dales, S. R. Patel and C. J. Read, Fr\'{e}chet algebras of power series, In Banach algebras 2009, Banach Center Publi. 91 (2010), 123-158.

\bibitem[12]{12} J. E. Forn\ae ss and N. \O vrelid, Finitely generated ideals in $A(\Omega)$, Ann. Inst. Fourier (Grenoble) 33 (1983), 77-85.

\bibitem[13]{13} O. Forster, Uniqueness of topology in Stein algebras, In Function algebras (ed F. T. Birtel) (Scott, Foresman, Glenview, Illinois, 1966), 157-163.

\bibitem[14]{14} A. M. Gleason, Finitely generated ideals in Banach algebras, J. Math. Mech. 13 (1964), 125-132.

\bibitem[15]{15} H. Goldmann, Uniform Fr\'{e}chet algebras, (North Holland Publ. Co., Amsterdam, 1990).

\bibitem[16]{16} N. Kerzman and A. Nagel, Finitely generated ideals in certain function algebras, J. Funct. Anal. 7 (1971), 212-215.

\bibitem[17]{17} I. Lieb, Die Cauchy-Riemannschen Differentialgleichung auf streng pseudokonveksen Gebieten: Stetige Randwerte, Math. Ann. 199 (1972), 241-256.

\bibitem[18]{18} R. J. Loy, Local analytic structure in certain commutative topological algebras, Bull. Austral. Math. Soc. 6 (1972), 161-167.

\bibitem[19]{19} R. J. Loy, Banach algebras of power series, J. Austral. Math. Soc. 17 (1974), 263-273.

\bibitem[20]{20} E. A. Michael, Locally multiplicatively convex topological algebras, Mem. Amer. Math. Soc. 11 (1952).

\bibitem[21]{21} J. B. Miller, Analytic structure and higher derivations on commutative Banach algebras, Aequationes Math. 9 (1973), 171-183.

\bibitem[22]{22} A. Noell, The Gleason problem for domains of finite type, Complex Variables Theory Appl. 4 (1985), 233-241.

\bibitem[23]{23} S. R. Patel, Fr\'{e}chet algebras, formal power series, and automatic continuity, Studia Math. 187 (2008), 125-136.

\bibitem[24]{24} S. R. Patel, Fr\'{e}chet algebras, formal power series, and analytic structure, J. Math. Anal. Appl. 394 (2012), 468-474.

\bibitem[25]{25} S. R. Patel, Uniqueness of the Fr\'{e}chet algebra topology on certain Fr\'{e}chet algebras, Studia Math. 234 (2016), 31-47.

\bibitem[26]{26} T. T. Read, The powers of maximal ideals in a Banach algebra and analytic structure, Trans. Amer. Math. Soc. 161 (1971), 235-248.

\bibitem[27]{27} S. J. Sidney, Properties of the sequence of closed powers of a maximal ideal in a sup-norm algebra, Trans. Amer. Math. Soc. 131 (1968), 128-148.

\bibitem[28]{28} M. P. Thomas, Local power series quotients of commutative Banach and Fr\'{e}chet algebras, Trans. Amer. Math. Soc. 355 (2003), 2139-2160.

\bibitem[29]{29} O. Zariski and P. Samuel, Commutative algebra. Vol. 2 (University Series in Higher Math., Van Nostrand, Princeton, New Jersey, 1960).
\end{thebibliography}
\end{document}